\documentclass[12pt]{article}
\usepackage{amsfonts}
\usepackage{amssymb}
\usepackage{amscd}

\def\Aff{\mathop{\mathrm {Aff}}\nolimits}
\def\aff{\mathop{\mathrm {aff}}\nolimits}
\def\Aut{\mathop{\mathrm {Aut}}\nolimits}
\def\Ad{\mathop{\mathrm {Ad}}\nolimits}
\def\ad{\mathop{\mathrm {ad}}\nolimits}
\def\Lie{\mathop{\mathrm {Lie}}\nolimits}
\def\PHC{\mathop{\mathrm {PHC}}\nolimits}
\def\HC{\mathop{\mathrm {HC}}\nolimits}

\def\HH{\mathop{\mathrm {HH}}\nolimits}
\def\Id{\mathop{\mathrm {Id}}\nolimits}
\def\Aff{\mathop{\mathrm {Aff}}\nolimits}
\def\aff{\mathop{\mathrm {aff}}\nolimits}
\def\Tot{\mathop{\mathrm {Tot}}\nolimits}
\def\SL{\mathop{\mathrm {SL}}\nolimits}
\def\sl{\mathop{\mathfrak {sl}}\nolimits}
\def\gl{\mathop{\mathfrak {gl}}\nolimits}
\def\SO{\mathop{\mathrm {SO}}\nolimits}

\def\tr{\mathop{\rm {Trace}}\nolimits}

\newtheorem{theorem}{Theorem}[section]
\newtheorem{proposition}[theorem]{Proposition}
\newtheorem{lemma}[theorem]{Lemma}
\newtheorem{remark}[theorem]{Remark}

\begin{document}
\title{K-theory and periodic cyclic homology  of some noncompact quantum algebras\footnote{Preliminary Version of July 10, 2003}}
\author{Do Ngoc Diep and Aderemi O. Kuku}
\maketitle
\begin{abstract} We prove in this paper that the periodic cyclic homology of the quantized algebras of functions on coadjoint orbits of connected and simply connected Lie group, are isomorphic to the periodic cyclic homology  of the quantized algebras of functions on coadjoint orbits of compact maximal subgroups, without localization. 
Some noncompact quantum groups and algebras were constructed and their irreducible representations were classified in recent works of Do Ngoc Diep and Nguyen Viet Hai \cite{diephai1}-\cite{diephai2} and Do Duc Hanh \cite{doduchanh} by using deformation quantization. In this paper we compute their K-groups, periodic cyclic homology groups and their Chern characters. 
\end{abstract}
\section*{Introduction} Let $G$ be a connected Lie group which acts smoothly on a locally convex algebra $A$ over the complex numbers, and $K$ a fixed maximal compact subgroup of $G$. 
One of the major results  of  V. Nistor in \cite{nistor1} is Theorem 1.1 , saying that up to localization at some maximal ideal $\mathfrak m$ of the algebra $C^\infty_{inv}(G)$ of bi-invariant functions $f$, (satifying $f(\gamma^{-1}x\gamma) \equiv f(x)$, $\forall g,x\in G$) the periodic cyclic homology of the crossed product $A \rtimes G$ and that for $A \rtimes K$, are isomorphic,  i.e.  if $q = \dim (G/K)$, then $$\PHC_*(A \rtimes G)_{\mathfrak m} \cong \PHC_{*+q}(A \rtimes K)_{\mathfrak m}. \eqno{(*)}$$ 

In the case where $A=\mathbb C$ and the action of $G$ on $\mathbb C$ is trivial, the crossed product $\mathbb C \rtimes G := C^\infty_c(G)$ becomes convolution algebra $C^\infty_c(G)$ and the restrictions of elements of $C^\infty_c(G)$, as functions on $G$, to the coadjoint orbits give  algebras of quantized functions on the coadjoint orbits. We prove the isomorphism (*) in that case without any localization. The main reason is that on coadjoint orbits, the bi-invariant functions correspond to functions which are constant on each coadjoint orbit and their localized forms are the same constants. (Note that  if we localize an algebra of functions  at a maximal ideal  $\frak P$, consisting of functions vanishing at a point, then an element of the algebra localized at the ideal  has the form of a formal series $$a = a_0 + a_1 + a_2 + \dots ,$$ where $a_0 \notin \mathfrak P$, $a_1 \in \mathfrak P$ but $a_1 \notin \mathfrak P^2$, $a_2 \in \mathfrak P^2$ but $ a_2 \notin \mathfrak P^3$, etc ..... and $a_0$ is uniquely defined (mod $\mathfrak P$), $a_1 = a - a_0,$ $a_2 = a - a_1$ (mod $\mathfrak P^2$), etc.... So, if a function is constant at a point then the localized form of that function at that point is exactly the constant value $a_0$.)
This idea will be put to use in the first section of this paper. It is interesting that with this isomorphism (without any localization) for the quantized algebras of functions on coadjoint orbits, we can reduce the computation of the noncommutative Chern-Connes characters of the Lie groups under consideration to those of maximal compact subgroups that are more easily computable.  Our main observation is that conjugacy in $G$ corresponds to the adjoint action on $\mathfrak g$ and to the coadjoint action on $\mathfrak g^*$. It is especially important in the concrete cases considered  in the last two sections, involving the Lie groups $\Aff(\mathbb R)$ and $\Aff(\mathbb C)$, (see \S3 and \S4) and the special linear group $SL(2,\mathbb R)$ in \S5.

\par
Indeed, the homogeneous classical mechanical systems with fixed Lie groups of symmetry were classified as coadjoint orbits of the Lie groups of symmetry or their central extensions by $\mathbb R$, in the vector space dual to the Lie algebras, see \cite{kirillov1}. For some special cases, where all the nontrivial orbits are of dimension equal to the dimension of the group (class $\overline{\rm MD}$), all such Lie groups have been classified and all the orbits explicitly computed. They are reduced to the cases of the groups of all affine transformations of the real or complex lines, see \cite{diep1}.  

The group of affine transformations of the real line has two 2-dimensional coadjoint orbits: the upper and lower half-planes, see \cite{diep1}. Using deformation quantization, a quantum analogue of the half-planes was constructed in \cite{diephai1}. The group of affine transformations of the complex line  has one orbit of complex dimension 2:  the punctured (withdraw a complex line through the origin) complex plane. Its quantum algebra was constructed in \cite{diephai2}. We compute the K-groups, the periodic cyclic homology of these quantum algebras and the corresponding Chern-Connes characters.
The groups $\SL(2,\mathbb R)$ has a family of elliptic hyperboloid orbits, a family of two-fold elliptic hyperboloid orbits, upper and lower elliptic half-cones, and the origin as a one-point orbit. The quantum algebras of functions over them are calculated in \cite{doduchanh}. In \S5 we also compute the K-theory, the periodic cyclic theory and the Chern-Connes characters.

In order to obtain these results, we use the methods from \cite{dkt1}-\cite{dkt2}. We first construct some diffeomorphims realizing canonical coordinates on coadjoint orbits and then reduce the quantum algebras to the ones related to the corresponding quantum groups. We then reduce these quantum algebras of quantum function on coadjoint orbits to the quantum algebras related to the maximal compact subgroups, that are more easily computable. This computation is realized in the last two sections.
\vskip .2cm
{\bf Notes on  Notation:} As usual we denote by capital letters some Lie groups, namely $G$, $K$, etc. Their corresponding Lie algebras will be denoted by the corresponding Gothic letters, namely $\mathfrak g$, $\mathfrak k$, etc. The dual space to a Lie algebra or a vector space will be denoted by the same letter with $*$, e.g. $\mathfrak g^*$ or $V^*$ is the dual space of the Lie algebra $\mathfrak g$ or the vector space $V$. $\PHC_*$ will denote the periodic cyclic homology, following \cite{nistor1}. If $A$ is a locally convex $\mathbb C$-algebra on which a Lie group $G$ acts smoothly, then $A \rtimes G$ denotes the crossed product of $A$ with $G$. $\Im(z)$ and $\Re(z)$ will mean the imaginary and real parts of the complex number $z$, while by $[a]$ and $\{a\}$  we denote the integral and fractional parts of $a$. $\mathbb R$ and $\mathbb C$ means the field of real or complex numbers, respectively. For any $\mathbb C$-algebra A, denote $A^\natural = \{A_n = A ^{\otimes n+1} \}$ the well-known Connes-Tsygan complex.

\section{Localization and coadjoint orbits}
Let $G$ be a connected and simply connected Lie group and $C^\infty_{inv}(G)$ the convolution algebra of bi-invariant functions on $G$, see \cite{nistor1}.
We prove in this section that localization of the convolution algebra $C^\infty_{inv}(G)$ at a maximal ideal corresponds to the (quantized) convolution algebra of functions with compact support on the corresponding orbit. We first modify some results obtained in the work \cite{nistor1} of V. Nistor.

\subsection{Preparation}
 Recall that a {\it quasi-cyclic object} in an Abelian category $\mathcal M$ is a graded object $(X_n)_{n \geq 0}$, $X_n \in Ob(\mathcal M)$ together with morphisms $d_i: X_n \to X_{n-1}$, for $i=0,\dots,n$ and $T_{n+1}: X_n \to X_n$ satisfying the following two axioms:
$$d_id_j =d_{j+1}d_i, \mbox{ for } i < j \leqno{(S1)}$$
$$d_iT_{n+1} = \left\{ \begin{array}{rl} T_n d_{i-1} &\mbox{ for } 1 \leq i \leq n\\   d_n &\mbox{for } i=0 \end{array}\right. \leqno{(C1)}$$
see \cite{nistor1} for more details. V. Nistor pointed out the examples of quasi-cyclic objects like: 
\begin{enumerate}
\item[(i)] the cyclic objects,
\item[(ii)] $\mathfrak A^\natural$, where $$\mathfrak A := A \rtimes G = C^\infty_c(G,A) = \{ \varphi \in C^\infty(G,A);  \varphi \mbox{ has compact support }\},$$ 
$G$ a Lie group,
$A$ a locally convex algebra, ($1 \in A$), on which $G$ {\it acts smoothly}, with the {\it twisted convolution product}
$$\varphi * \psi (g) := \int_G \varphi(h) \alpha_h(\psi(h^{-1}g)dh,$$ 
where $dh$ is a fixed left invariant Haar measure on $G$, $\alpha: G \to \Aut(A)$ a {\it smooth action} of $G$ on $A$ in the sense that the map $g \mapsto \alpha_g$ is continuous and unital and the map $g \mapsto \alpha_g(a)$, for all $a\in A$ is smooth. In that case, we have $\mathfrak A^\natural_n = (A \rtimes G)^\natural_n = C^\infty_c(G^{n+1}, A^{\otimes n+1})$ with the operations
$$(d_j\varphi)(g_0, g_1, \dots, g_{n-1}) :=$$ $$\int_G d_j \circ (1 \otimes \dots \otimes 1 \otimes \alpha_h \otimes 1 \otimes \dots \otimes 1)(\varphi(g_0, \dots, g_{j-1}, h, h^{-1}g_j, \dots, g_{n-1}))dh,$$ $j=0,1,\dots, n-1$,
$$(d_n\varphi)(g_0, g_1, \dots, g_{n-1}) :=$$ $$\int_G d_n \circ (\alpha_h \otimes 1 \dots \otimes 1)(\varphi(h^{-1}g_0, g_1,\dots,  g_{n-1},h))dh$$
and 
$$(t_{n+1}\varphi)(g_0,\dots,g_{n}) := t_{n+1}(\varphi(g_1,\dots,g_{n},g_0))$$
As remarked in \cite{nistor1}, the operators $d_j, t_{n+1}$ on the right hand side are the ones of the cyclic Connes-Tsygan complex $A^\natural$ (see \cite{connes}).

In particular, if $A= \mathbb C$ and the action of $G$ on $\mathbb C$ is trivial, we have the convolution algebra $C^\infty_c(G)$. For the applications in this paper, our attention will be forcused on $C^\infty_c(G)$.
\item[(iii)] Let $A$, $G$ be as in (ii) above. $G_1$ some other Lie group, $U$ an open subgroup of $G$. Write $L_n(U,G_1) := C^\infty_c(U \times G_1^{n+1}, A^{\otimes n+1})$. Then $\{L_n(U,G_1)$ is a quasi-cyclic object with similar operations $d_j, T_{n+1}$, (see \cite{nistor1} for details).
\end{enumerate}

\subsection{An $H$-relative cohomology complex}
Let us now introduce some new examples of quasi-cyclic objects, related with some quotient maps. In (ii) above, $G_1$ is a group, $\rho$ a homomorphism. However in what follows, $G_1$ is replaced by the homogeneous space $H\setminus G$ and $\rho$ is just the quotient map. 

Consider the canonical quotient map $\rho: G \to H\setminus G$, for some subgroup $H$. Consider an open set $U \subset G$. Define 
$$L_n(U,H\setminus G) := C_c^\infty(U \times (H\setminus G)^{n+1}, A^{\otimes n+1})$$ and define also 
$$(d_j\varphi)(\gamma, H g_0,\dots, \hat{g_j}, \dots,H g_n) :=\int_G d_j (1 \otimes \dots \otimes \alpha_\gamma \otimes \dots \otimes 1) (\varphi(\gamma, Hg_0, \dots, Hg_{n}))d\mu(g_j),$$ $j=1,\dots, n$, $d\mu(.)$ is the quotient measure on the quotient space $H\setminus G$, 
$$(d_0\varphi)(\gamma, \hat{g_0}, Hg_1, \dots, Hg_n) :=\int_G d_0 (1 \otimes \alpha_\gamma \otimes \dots \otimes 1)(\varphi(\gamma, Hg_0, Hg_1,\dots,  Hg_{n}))d\mu(g_0)$$
and 
$$(T_{n+1}\varphi)(\gamma, Hg_0,\dots,Hg_{n}) := (1\otimes\alpha_\gamma^{-1} \otimes \dots, \otimes 1)t_{n+1}(\varphi(\gamma, Hg_1,\dots,Hg_{n},Hg_0))$$
 
\begin{proposition}
$(\{L_n(U,H\setminus G)\}, d_j, d_0, T_{n+1})$ is a quasi-cyclic object.
\end{proposition}
{\it Proof}.
By similar arguments to those in the work of V. Nistor, see \cite{nistor1}, it is easy to see that we have also a quasi-cyclic object. \hfill$\Box$

This quasi-cyclic object is related to the quantum algebras of functions on orbits, as we shall see later.
\vskip .2cm
 
Let us now define an action of $G$ on the quasi-cyclic object $(\{L_n(U,H\setminus G)\}, d_j, d_0, T_{n+1})$. For a fixed action $\alpha : G \to \Aut(A)$, define $\beta: G \to GL(L_n(U,H\setminus G))$ by
$$\beta_\gamma(\varphi)(\gamma, Hg_0, Hg_1, \dots, Hg_n) := \alpha_\gamma^{\otimes n+1}(\varphi(\gamma^{-1}\gamma_1\gamma, Hg_0\gamma, \dots, Hg_n\gamma)),$$ for all $\gamma$ in $G$ and $ \gamma_1$ in $U$, $g_0, \dots, g_n$ in $G$, $\varphi$ in $L_n(U,H\setminus G)$.

In particular, if $H=\{ e \}$ and $U=G$, we have $H\setminus G \cong G$ and a map $p: L_n(G,G) \to (A \rtimes G)^\natural_n$, defined by
$$p\varphi(h_0, \dots, h_n) := \int_G \Psi(\beta_\gamma\varphi)(g_n,g_0,\dots,g_{n-1}),$$ where $g_0 := h_0$, $g_1 := h_0h_1$, $g_2 :=h_0h_1h_2$, .... $g_n := h_0h_1h_2\dots h_n$, $\Psi := (\alpha_{g_n} \otimes \alpha_{g_0} \otimes \dots \alpha_{g_{n-1}})^{-1}$.

\vskip .2cm
We now define the map $J$ which gives rise to an isomorphism in Hochschild homology, (see \cite{nistor1}).  

Consider again an open set $U$ which is $\Ad_G$-invariant in $G$, and define $J: L_n(U, H\setminus G) \to L(U, \{e\})$ by 
$$(J\varphi)(\gamma) := \int_{(H\setminus G)^{n+1}} \varphi(\gamma, Hg_0,\dots, Hg_n)d\mu(g_0)\dots d\mu(g_n)$$
Note that we use $d\mu(g)$  here to denote the relative quasi-invariant measure on the quotient space $H\setminus G$.
 
\begin{lemma} $J$ is a morphism of quasi-cyclic objects and $HH(J)$ is an isomorphism of the corresponding Hochschild homology groups , i.e. $\HH_n(L(U,H\setminus G) \cong \HH_n(L(U,\{e\})$.\end{lemma}
{\it Proof.} As in \cite{nistor1}, since we have only replace $G_1$ by $H\setminus G$ in \cite{nistor1}.
\hfill$\Box$
\vskip .2cm

Suppose we have some continuous homomorphism of a compact group $K$ into $G$. We have then 
\begin{proposition}
$\HC_*(L(K,H\setminus G)^K) \cong \HC_*(((A \rtimes K)^\natural)^H)$
\end{proposition}
{\it Proof}. We have from \cite{nistor1} the isomorphism
$\HC_*(L(K,G)^K) \cong \HC_*(((A \rtimes K)^\natural))$. From the definition of complexes defining $\HC_*(L(K,G)^K)$ and $ \HC_*(((A \rtimes K)^\natural))$, we have isomorphic $H$-invariant homology groups $\HC_*((L(K, G)^K)^H) \cong \HC_*(((A \rtimes K)^\natural)^H)$. But $\HC_*((L(K, G)^\natural)^H) \cong \HC_*(L(K,H\setminus G)^K)$
\hfill$\Box$

\vskip .2cm
Let $G$ be a connected Lie group, $\mathcal F$ a smooth $G$-module. Let $K\subset G$ be a maximal compact subgroup of $G$, $\mathfrak k = \Lie K$, $\mathfrak g = \Lie G$ the corresponding Lie algebras,  $q=\dim(G/K)$. 
We now define $G$-equivariant homology: Consider the complex of relative Lie algebra homology 
$$0 \to (\wedge^q(\mathfrak g/\mathfrak k) \otimes_H \mathcal F) \otimes_K \mathbb C {\buildrel\delta \over \longrightarrow} \dots {\buildrel\delta \over \longrightarrow} (\wedge^0(\mathfrak g/\mathfrak k) \otimes_H  \mathcal F) \otimes_K \mathbb C \longrightarrow  \mathcal F \otimes_G \mathbb C \longrightarrow 0,\eqno{(I)}$$
 where $\delta: (\wedge^j(\mathfrak g/\mathfrak k) \otimes_H \mathcal  F) \otimes_K \mathbb C \longrightarrow (\wedge^{j-1}(\mathfrak g/\mathfrak k) \otimes_H \mathcal F) \otimes_K \mathbb C$ is defined as
$$\delta(\dot X_1 \wedge \dots \wedge \dot X_j \otimes \xi) := \sum_{i=1}^j (-1)^{i+1} \dot X_1 \wedge \dots \wedge \hat{X_i} \wedge \dots \dot X_j \otimes X_i(\xi) -$$ $$- \sum_{i < j} (-1)^{i+k}[X_i,X_k] \wedge \dot X_1 \wedge \dots \wedge \hat{X_i} \wedge \dots \wedge \hat{X_k}\wedge\dots\wedge \dot X_j \otimes \xi,$$
and for $X_i\in \mathfrak g, \dot X_i$ is the class of $X_i$ in  $\mathfrak g/\mathfrak k$, and $\hat X_i$ indicates that  $\dot X_i$ is omitted. In the complex (I) $\mathcal F$ is regarded as a smooth $H$-module.
It is not hard to prove
\begin{proposition}
The complex (I) is {\it acyclic}.
\end{proposition}

For acyclicity of a similar complex with $\mathbb C$ instead of  $H$, see \cite{nistor1}(Prop. 3.6). In our relative case this is true also, as can be verified.

\subsection{A relative complex}
In the main body of the paper of V. Nistor is the construction of complexes computing the mentioned isomorphism,   $$\PHC_*(A \rtimes G)_\mathfrak m \cong \PHC_{*+q}(A \rtimes K)_\mathfrak m.$$
In this section, we shall introduce a relative complex satisfying all the conditions of the Main Lemma of Nistor (see 1.3.2, 1.3.3).  First we briefly review the conditions satisfied by Nistor's lemma to fix notations (see 1.3.1).

\subsubsection{}
Following \cite{nistor1}, let us consider the following {\bf data}:
\begin{itemize}
\item Two exact sequences of quasi-cyclic objects of complete locally convex spaces
$$\leqno{(E1)}\qquad\CD 0@<<<\mathcal Y @<<< \mathcal{X}^{(0)} @<\delta<< \dots @<\delta<< \mathcal X^{(q-1)}@<\delta<< \\ @. @. @.  @<\delta<< \mathcal X^{(q)} @<<< 0 \endCD$$
$$\leqno{(E2)}\qquad\CD 0 @>>> \mathcal{X}^{(0)} @>\sigma>> \dots @>\sigma>> \mathcal X^{(q-1)} @>\sigma>> \mathcal X^{(q)} @>>> \\ @. @. @. @>>> \mathcal Z@>>> 0 \endCD$$
\item A $C^\infty$-action of $\mathbb R$ on $\mathcal X^{(j)}$ for any $j$ such that $\eta_1 = T^{n+1}_{n+1}$ and the derivative $\nabla = \frac{d\eta_t}{dt}\vert_{t=0}$ of $\eta$ at $t=0$ is equal to $\nabla=\delta\sigma +\sigma\delta$, with the convention that $\delta(\mathcal X^{(0)})= \sigma(\mathcal X^{(q)}) = 0$
\item The endomorphisms $1-T_{n+1}^{n+1}$ are injective on $\mathcal X^{(j)}_n$ for all $j=0,\dots,q$ and for all $n \geq 0$. 
\end{itemize}

Let us recall that a {\it precyclic object} is a quasi-cyclic object $((X_n)_{n \geq 0}, d_j, T_{n+1})$ such that $T_{n+1}^{n+1} = 1$. Given a precyclic object $X$, one constructs the {\it Connes-Tsygan complex} $\mathcal C(X)$ as in the cyclic case, see \cite{nistor1}. The homology of the {\it 2-periodic total complex} $\Tot\mathcal C(X)$ associated to the bi-complex $\mathcal C(X)$ is defined as the periodic cyclic homology $\PHC(X)$ of the complex $X$, (see \cite{nistor1}, Definition 2.2, and the definition thereafter).  

We now associate to the data satisfying the above definition, some new objects 
$$\tilde{\mathcal X}^{(0)} := \mathcal X^{(0)}/\nabla\mathcal X^{(0)}, \quad
\tilde{\mathcal X}^{(j)}:= \mathcal X^{(j)}/(\nabla \mathcal X^{(j)} + \sigma \mathcal X^{(j-1)}), \forall j=1,\dots , q. $$

and state the followings Lemma due to V. Nistor (see \cite{nistor1} 
\begin{lemma}
\begin{enumerate}
\item[(i)] For any $j=1,\dots, q$, the object $\tilde{\mathcal X}^{(j)}$ is a precyclic object.
\item[(ii)] The complex
$$\CD 0 @<<< \mathcal Y @<<< \tilde{\mathcal X}^{(0)} @<\delta<< \tilde{\mathcal X}^{(1)} @<\delta<<\dots \\ @. \dots
@<\delta<< \tilde{\mathcal X}^{(q-1)} @<\delta<< \mathcal Z \cong \mathcal X^{(q)}/\sigma\mathcal X^{(q-1)} @<<< 0 \endCD$$ is acyclic.
\item[(iii)] $\PHC_*(\tilde{\mathcal X}^{(j)}) = 0$ for any $j=0,\dots,q$.
\item[(iv)] $\PHC_*(\mathcal Y) \cong \PHC_{*+q}(\mathcal Z).$
\end{enumerate}
\end{lemma}
We now verify conditions of this lemma for the $H$-relative complex.
 
\subsubsection{Relative form of $\mathcal{X}^{(j)}$ and $\delta$}
For any Lie group $G$, let $\mathbb C_\Delta$ be the one-dimensional representation of $G$ by multiplication with the modular function $\Delta$, $G_x = \{ \gamma \in G \vert \gamma x = x\gamma \}$, $\mathfrak g_x = \Lie G_x$, $K$ a maximal compact subgroup of $G$ and $U \subset G$, $K_x = K \cap G_x$.
Let $\Delta'$ denote the modular function for $G_x$, $U'$ an open subset of $G_x$,
$$\mathcal F' := \left\{ \begin{array}{ll}
\mathbb C_{\Delta'} \otimes L(U',G_x), &\mbox{if } x \in K_x,\\
\mathbb C_{\Delta'} \otimes L(U',G_x)/ (1-x)(C_{\Delta'} \otimes L(U',G_x)), &\mbox{if } x \notin K_x\end{array}\right. $$
 $$M := \left\{\begin{array}{ll}
K_{x}, & \mbox{if } x \in K_x,\\
\mbox{The maximal compact subgroup in }G_{x}/(x),& \mbox{if } x \notin K_x
\end{array}\right. $$ and 
$\mathfrak m := \Lie(M)$.

Define the $H$-relative ${\mathfrak g}_x$-cohomology complex $\{C_j^H\}$ of $\mathcal F'$ 
by $C_j^H := (\wedge^j {\mathfrak g}_x \otimes \mathcal F')^H$ and $\delta_0: C_j \to C_{j-1}$ with
$$\delta_0(X_0 \wedge\dots\wedge X_j \otimes \xi) = \sum_{j=1}^j (-1)^{i+1} X_1 \wedge\dots\wedge\hat{X}_i\wedge \dots\wedge X_j \otimes X_i(\xi)$$
$$- \sum_{1\leq i < k \leq j}(-1)^k[X_i,X_k]\wedge X_1\wedge\dots\wedge\hat{X_i}\wedge\dots\wedge X^k\wedge\dots\wedge X_j \otimes \xi$$
Now we define $\mathcal X^{(j)} := C_j/C'_j$, where $C'_j$ is generated by $\mathfrak m \wedge C_{j-1} \mbox{ and } (\gamma-1)C_j$ for all $\gamma\in K_x \cup (x)$. It is not hard to see that 
$$\delta_0(X \wedge \omega) = X(\omega) - X\wedge \delta_0(\omega), \forall X \in \mathfrak g_x, \forall \omega \in C_{j-1},$$ where $X(\omega)$ means the contraction of $X$ and $\omega$ with values in $C_{j-1}$.  
Define $\delta$ to be the quotient map of $\delta_0$ on $\mathcal X^{(j)} := C_j/C_j'$. 

Define also $d_j: C_k^H \to C_k^H, j=1,\dots,k$ by
$$d_j(X_1\wedge\dots\wedge X_k \otimes \lambda\otimes\xi) := X_1\wedge\dots\wedge X_k \otimes \lambda \otimes d_j\xi, $$ $\forall X_1\wedge\dots\wedge X_k\in \wedge^k \mathfrak g_x, \forall \lambda\in \mathbb C_\Delta, \forall \xi\in \mathcal F$.

\subsubsection{Definition of $\sigma$ in the relative case}
Let $C_j$ $\sigma$ be as in 1.3.1 and 1.3.2. 
For every element $Z\in C^\infty(U',\mathfrak g_x)$, define $Z(x\exp X) := X$ and $$\sigma_0: C_k \to C_{k+1}; \sigma_0(\omega) := Z \wedge \omega,$$
Define $\sigma : \mathcal X^{(k)} \to \mathcal X^{(k+1)}$ to be the quotient maps of $\sigma_0$ mode $C'_j$.

\subsubsection{Relative form of $\eta'$ and the equation $\nabla = \sigma\delta + \delta\sigma$}
Let $G_x$, $H$ be as before, $H' = H \cap G_x$.
First define an action $\eta'$ of $\mathbb R$ on $L(U',H'\setminus G_x)$ as 
$$(\eta'_t\varphi)(x\exp X, H' g_0, \dots, H' g_n) := \beta_{\exp(tX)}(x\exp(X),H' g_0, \dots, H' g_n),\eqno{(**)}$$ for any $\varphi\in L(U',H'\setminus G_x)$, $t\in \mathbb R$. It is natural to extend this action to $C_j$ by 
$$t.(X_1\wedge\dots\wedge X_k\otimes\lambda\otimes\xi) := X_1\wedge\dots\wedge X_k\otimes \lambda \otimes \eta'_t(\xi),$$
for all $X_1,\dots,X_k\in \mathfrak g_x$, $\lambda\in \mathbb C_\Delta$ and $\xi\in L(U',H\setminus G_x)$. Let $\nabla_0$ be the derivative of the action of $\eta_t'$ above at 0. By a similar computation as in
\cite{nistor1}, one can see that  
$$\nabla_0 = \delta_0\sigma_0 + \sigma_0\delta_0.$$  
It is also not hard to see that each $C'_j$ is invariant under this action and therefore we get and action $\eta$ of $\mathbb R$ in $\mathcal X^{(j)}$.

Finally for the above data, in a way analogous to that in \cite{nistor1}, we define 
$${\mathcal Y} := (\mathbb C_{\Delta'} \otimes_H L(U',H'\setminus G_x)) \otimes_{G_x} \mathbb C = ((A \rtimes G_x)^\natural_{U'})^H$$
and 
$$\mathcal Z := \left\{\begin{array}{ll}
L(U' \cap K_x, H'\setminus G_{x}) \otimes_{K_{x}} \mathbb C & \mbox{if } x\in K_x\\
L(U' \cap M_0, H'\setminus G_{x}) \otimes_M \mathbb C & \mbox{if } x \notin K_x
\end{array}\right.,$$ where $M_0$ is the inverse image of $M$ in $G_{x}$ of the maximal compact subgroup of $G_{x}/(x)$.

By a similar argument to that in \cite{nistor1}, we can conclude that  for the $H$-relative complex all the conditions of the main lemma of \cite{nistor1} are also satisfied. 

\subsection{Passage to coadjoint orbits}
After the statement of his main result (Theorem 1.1, p. 4 in \cite{nistor1}), V. Nistor stated as follows. ``This fits with Mackey's method of orbits, except that now for reasons we do not yet understand, we obtain orbits on $(\Lie G)$ rather than in $(\Lie G)^*$.  An interesting feature of the result is worthwhile stressing: there is no $\gamma$-obstruction in cyclic cohomology.'' This means that he didn't work with the coadjoint orbits. 
We now explain that it is natural to pass to coadjoint orbits and that the localization disappears on coadjoint orbits. 

Let $G$ be a connected and simply connected Lie group, $\mathfrak g = \Lie(G)$, $\mathfrak g^*$ the dual of $\mathfrak g$. The Lie algebra, by the well-known Ado Theorem {\it can be considered as linear} , i.e. as a subalgebra of some general linear algebra $\gl(n,\mathbb C)$ and therefore the $\Ad$ action becomes conjugation, i.e. $\mathfrak g \hookrightarrow \gl(n,\mathbb C)$ and for $X\in \mathfrak g$, $g\in G$, $\Ad_G(g)X = gXg^{-1}$. Let $F$ be a fixed point in $\mathfrak g^*$, $G_F$ the stabilizer of $F$.
Let us consider the natural projection $G \twoheadrightarrow G_F \setminus G = \Omega $ defined as $x_0 \mapsto G_Fx_0 $ where $G_Fx_0$ corresponds to $ \tilde{F}_0 \in \mathfrak g$ inder a fixed linear isomorphism $\mathfrak g \cong \mathfrak g^*$, let $F_0\in \mathfrak g^*$ corresponds to $\tilde{F_0}$ in $\mathfrak g$. Suppose that $x_0 = \exp(\tilde{F_0})$, i.e. $\tilde{F_0} = \ln x_0$, then from the well-known Van Campbell-Hausdorff-Dynkin formula for $\ln(\exp(X)\exp(Y))$, we deduce that $\ln(xx_0x^{-1}) = \Ad_x\tilde{F_0}$. We have therefore the following
\begin{lemma} Under the map $G \twoheadrightarrow G_F \setminus G = \Omega \hookrightarrow  \mathfrak g^*$, the element $xx_0x^{-1}$ in $G$ goes to the element $\Ad(x^{-1})F_0$, and the conjugacy orbit of $x_0$ goes to the coadjoint orbit 
$$\Omega_{F_0} = \{ \Ad(x^{-1})F_0 \quad\vert\quad x\in G \}.$$
\end{lemma}

\subsection{Localized functions on coadjoint orbits}
We apply the construction of the subsection 1.2 to the case of coadjoint orbit $\Omega = G_F \setminus G$.
\begin{lemma} \label{lemma1.1} 
Let $G$ be a connected and simply connected Lie group. There is a vector space isomorphism $C^\infty_c(\mathfrak g) \cong C^\infty_c(U)$, where $U$ is an open set in $G$.
\end{lemma} 
{\it Proof.} The exponential map $\exp: \mathfrak g \to G$ is a local diffeomorphism and  we can choose an open neighborhood $V$ of $X$ in $\mathfrak g$ and an open neighborhood $U=\exp(V)$ of $x=\exp(X)$ in $G$ such that $\exp : V \to U$ is a diffeomorphism. We have a map $C^\infty_c(U) \to C^\infty_c(\mathfrak g)$. If some function $f\in C^\infty_c(U)$ is such that $f(\exp X) = 0$ for all $X\in \mathfrak g$, then because it has compact support, we can find a finite  covering of its support by the open neighborhood of above type where the exponential map becomes diffeomorphism. It is therefore identically $0$ and the map  $C^\infty_c(U) \to C^\infty_c(\mathfrak g)$ is an an injection.

For an arbitrary element $f\in C^\infty_c(\mathfrak g)$ we can cover its support by a finite number of the above neighborhoods where the exponential map became diffeomorphic. The composition of $f$ with $\exp^{-1}$ should be a function $\varphi$ on $U$ such that $f = \varphi\circ \exp$. The map is therefore a surjection. 
\hfill$\Box$. 

\begin{lemma}\label{lemma1.2}
There is a isomorphism between convolution algebras 
$$ C^\infty_c(\mathfrak g) \cong C^\infty_c(\mathfrak g^*)$$
\end{lemma}

{\it Proof.} We can identify $\mathfrak g$ with $\mathfrak g^*$ by fixing some basis in $\mathfrak g$ and the corresponding dual basis in $\mathfrak g^*$. On the vector space $C^\infty_c(\mathfrak g^*)$ there is a  $\star$-product of functions, which is related with the coadjoint action of the Lie group, i.e. for any $f,g\in C^\infty_c(\mathfrak g)$ define the Possion structure $\{f,g\}$, as follows. For $X\in \mathfrak g^*$, $df(X), dg(X) \in \mathfrak g$ put 
$$\{f,g\}(X):= \langle X, [df(X),dg(X)]\rangle.$$ With this Poisson bracket,  is associated a natural $\star$-product of functions, see e.g. \cite{fedosov}.
We define the $\star$-product structure on $C^\infty_c(\mathfrak g)$ by transporting the $\star$-product structure on $C^\infty_c(\mathfrak g^*)$ through the fixed isomorphism $\mathfrak g \cong \mathfrak g^*$. 
\hfill$\Box$

\begin{lemma}\label{lemma1.3}
There is a isomorphism between convolution algebras 
$$ C^\infty_c(\mathfrak g)^G \cong C^\infty_c(\mathfrak g^*)^G$$
\end{lemma}
{\it Proof}. Under the isomorphism $\mathfrak g \cong \mathfrak g^*$, the adjoint action $\Ad_Gx$ becomes the coadjoint action $\Ad_G^*(x^{-1})$. Because we have an isomorphism of linear structures $ C^\infty_c(\mathfrak g)^G \cong C^\infty_c(\mathfrak g^*)^G$ and that there is a  $\star$-product on $C^\infty_c(\mathfrak g^*)$, we can define a $\star$-product on $C^\infty_c(\mathfrak g)$ from $C^\infty_c(\mathfrak g^*)$, using this isomorphism. Hence, we have an isomorphism $ C^\infty_c(\mathfrak g)^G \cong C^\infty_c(\mathfrak g^*)^G$.
\hfill$\Box$

\vskip .5cm
We turm our attention  from now on to linear Lie groups. For linear Lie groups, the Ad action become conjugation and the adjoint orbits become conjugacy classes.

\begin{remark}
1. For the connected and simply connected Lie groups of type I, there is one-to-one correspondecne between the irreducible unitary representations and the orbit data \cite{auslander_kostant}, see also \cite{kirillov1}, \cite{diep1}.

2. For all almost algebraique Lie groups Duflo had shown also that there are one-to-one correspondences between a big enough set (of Plancherel measure 0) of the irreducible unitary representations of $G$ and the orbit data and the central characters
\cite{duflo}.
\end{remark}

\begin{lemma}\label{lemma1.4}
There is a one-to-one correspondence between the localization of the algebra $C^\infty_{inv}(G)$ at maximal ideals and the central characters of the Lie group $G$, i.e. a bi-invariant function in $C^\infty_{inv}(G)$, localized at any point on a fixed coadjoint orbit $\Omega$ takes the same constant value as some central character of an irreducible unitary representation $T_{F,\sigma}$ corresponding to some orbit datum $(F,\sigma)$ associated with the orbit $\Omega$.  
\end{lemma}
{\it Proof.} 
Let us first recall what  ``localization'' means. For a fixed maximal ideal $\mathfrak P$ in an algebra $A$, an element $a\in A$ is localized at the point $\mathfrak P$ means a decompositon of $a$ into a series $$a = a_0 + a_1 + \dots $$ such that $ a_0 \notin \mathfrak P, a_1 \in \mathfrak P$ but $a_1 \notin \mathfrak P^2$, $a_2 \in \mathfrak P^2$ but $a_2 \notin \mathfrak P^3$, etc.... 
On the coadjoint orbits, bi-invariant functions $\varphi \in C^\infty_{inv}(G)$ are constant. Therefore, the value of a bi-invariant function, localized at any maximal ideal, (consisting  of all the functions vanishing at the corresponding point of the coadjoint orbit), is the same as the constant the bi-invariant function takes on the coadjoint orbit. The value of the constant function on coadjoint orbits are the same as the value of a the central character of the representation $T_{F,\sigma}$ associated to these orbit, (see for example \cite{kirillov1}, \cite{duflo}).\hfill$\Box$

\vskip .2cm
To be able to state the main theorem \ref{theorem1.7} below,
let us first fix some notations. For any Lie group $G$ and subgroups $K$, $\Omega = \Omega_F$ a coadjoint orbit of $G$, passing through a fixed point $F$ in $\mathfrak g^*$, we shall write $\Omega\vert_K = \Omega_{F\vert K}$ for the coadjoint K-orbit passing through $F\vert_{\mathfrak k}$. We also write $K_{F\vert_K}$ for the stabilizer of coadjoint $K$-action at $F\vert_{\mathfrak k}$.

\begin{theorem} \label{theorem1.7} Let $G$ be a connected and simply connected Lie group such that there are one-to-one correspodences between the central characters, the irreducible unitary representations and the coadjoint orbit data. Let $\Omega = \Omega_F$ a coadjoint orbit of $G$ passing through a fixed point $F$ in $\mathfrak g^*$, $K$ a maximal compact subgroup of $G$, $q=\dim(G/K)$ and $\Omega_K$ the coadjoint orbit passing through $F\vert_{\mathfrak k}$. Let $C^\infty_c(\Omega)$ (resp. $C^\infty_c(\Omega\vert_K))$ be the quantized algebra of functions on $\Omega$ (resp., $\Omega\vert_K$). Then 
we have an isomorphism
$$\PHC_*(C^\infty_c(\Omega_G) \cong \PHC_{*+q}(C^\infty_c(\Omega_K)).$$
\end{theorem}

{\it Proof.} First observe that  a coadjoint orbit $\Omega_F \subseteq \mathfrak g^*$ can be identified with the homogeneous space $H\setminus G$, where $H=G_F$ is the stabilizer of an arbitrary point $F$ in $\Omega_F$. We have therefore $$C^\infty_c(\Omega_F) \cong C^\infty_c(G_F\setminus G)\cong C^\infty_c(G)^H.\eqno{\star}$$ 
Let us write some object with sub-index $K$ to mean its restriction onto $K$. Also we have the same notations for a fixed maximal compact subgroup $K$, $C^\infty_c(\Omega_{F|_K}) \cong C^\infty_c(K_{F_K}\setminus K)\cong C^\infty_c(K)^{H_K}$. From Proposition 1.3, we have $$\HC_*(L(K,H\setminus G)^K) \cong \HC_*((A \rtimes K)^K)^H.$$ Also we have from Proposition 1.4 that the complex $$0 \to (\wedge^q(\mathfrak g/\mathfrak k) \otimes_H \mathcal F) \otimes_K \mathbb C {\buildrel\delta \over \longrightarrow} \dots {\buildrel\delta \over \longrightarrow} (\wedge^0(\mathfrak g/\mathfrak k) \otimes_H \mathcal F) \otimes_K \mathbb C \longrightarrow \mathcal F \otimes_G \mathbb C \longrightarrow 0$$ is acyclic. Using Lemma 1.9, we can identify $C^\infty_c(\Omega_F)$ with some algebra of functions on a conjugacy orbit of $G$ in $\mathfrak g$ with adjoint action. A function in $C^\infty_{inv}(G)$ means a function which is constant on conjugacy classes of $G$. So from Lemmas 1.7-1.9, it is the same as a function, which is constant on each coadjoint orbit. From Lemmas 1.11, a bi-invariant function localized at a maximal ideal $\mathfrak m$ in the space of central functions, takes the same value as some central character of some irreducible unitary representation, corresponding  to the coadjoint orbit. The localization at the ideal $\mathfrak m$ therefore means taking restrictions of functions on the associated orbit. Now we can apply the theorem 1.1 of Nistor\cite{nistor1} to conclude that  the homology groups with localization are isomorphic to the groups without localization.
\hfill$\Box$

\section{Noncommutative Chern-Connes characters}

Let us now briefly recall the construction of noncommutative Chern-Connes characters. We then present an interesting result (Theorem 2.1) for the noncommutative Chern-Connes characters of the quantized algebras of functions on coadjoint orbits.

\subsection{K-groups. Connes-Kasparov-Rosenberg Theorem}.
In the following, $K$-groups shall mean the $\mathbb Z/(2)$-graded algebraic $K$-groups of algebras over the field of complex numbers. For connected solvable Lie groups $G$ and  a maximal compact subgroup $K$, the so called Connes-Kasparov conjecture  has been proved (see, e.g. \cite{judith}) and it is known for large classes of Lie groups as Connes-Kasparov-Rosenberg Theorem, i.e.
$$K_*(C^\infty_c(G)) \cong K_{*+q}(C^\infty(K)),$$ where $q = \dim(G/K)$. 

\subsection{Noncommutative Chern-Connes characters}
For the general notion of Chern-Connes characters, readers are referred to the work of J. Cuntz \cite{cuntz}.  
From the works of J. Cuntz \cite{cuntz} and V. Nistor \cite{nistor1}-\cite{nistor2} we can deduce that there is a natural noncommutative Chern-Connes character $ch$ with values in the localizations of $\PHC_*(A)$, as $$ch: K_*(C^\infty_c(G)) \to \PHC_*(C^\infty_c(G))_{\mathfrak m}.$$ 

\subsection{The commutative diagrams}

\begin{theorem} Let $G$ be a connected and simply connected Lie group such that there are one-to-one correspodences between the central characters, the irreducible unitary representations and the coadjoint orbit data, $\Omega = \Omega_F$ a coadjoint orbit of $G$ passing through a fixed point $F$, $K$ a maximal compact subgroup of $G$, $\Omega\vert_K$ the coadjoint orbit of $K$ passing through $F\vert_{\mathfrak k}$, $\mathbb T$ the maximal torus of $G$ in $K$ and $W=N(\mathbb T)/\mathbb T$ the Weyl group corresponding to $\mathbb T$. For any co-adjoint orbit $\Omega$, let $C^\infty_c(\Omega)$ be the quantized algebra of functions  on $\Omega$ with compact support. Then, there 
 is a commutative diagram for the noncommutative Chern-Connes characters of the quantized algebra of functions on coadjoint orbits
$$\CD
K_*(C^\infty_c(\Omega)) @>ch_\Omega>> \PHC_*(C^\infty_c(\Omega)))\\
@VIVV   @VVIIIV\\
K_{*+q}(C^\infty_c(\Omega_K)) @>ch_{\Omega_K}>> \PHC_{*+q}(C^\infty_c(\Omega_K))\\
@VIIVV   @VV{IV}V\\
K_{*+q}(C^\infty_c(\mathbb T))^W @>ch_{\mathbb T}>> \PHC_{*+q}(C^\infty_c(\mathbb T))^W
\endCD$$
and modulo torsion, the noncommutative Chern-Connes characters are isomorphisms.
\end{theorem}
{\it Proof.} 
For a Lie group such that there are one-to-one correspodences between the central characters, the irreducible unitary representations and the coadjoint orbit data $G$ and a maximal compact subgroup $K$, there is a natural commutative diagram: 
$$\CD
K_*(C^\infty_c(G)) @>ch_G>> \PHC_*(C^\infty_c(G)))_{\mathfrak m}\\
@VVV   @VVV\\
K_{*+q}(C^\infty_c(K)) @>ch_K>> \PHC_{*+q}(C^\infty_c(K))_{\mathfrak m}\endCD$$
where the first vertical row is the Connes-Kasparov-Rosenberg isomorphism, and the second vertical one is the isomorphism of V. Nistor \cite{nistor1}. 
This can be reduced to the maximal torus case 
$$\CD
K_{*+q}(C^\infty_c(K)) @>ch_K>> \PHC_{*+q}(C^\infty_c(K))_{\mathfrak m}\\
@VVV   @VVV\\
K_{*+q}(C^\infty_c(\mathbb T))^W @>ch_{\mathbb T}>> \PHC_{*+q}(C^\infty_c(\mathbb T))^W_{\mathfrak m}\endCD$$ where $W = W(\mathbb T)$ is the Weyl group corresponding to the maximal torus $\mathbb T$ and the sub-indices of $ch$ indicate the corresponding target groups. 
In this second commutative diagram the first vertical row is the well-known result of the K-theory of compact Lie groups and the second one exists via Morita equivalence (see \cite{dkt1},\cite{dkt2}). The horizontal row at the bottom is an isomorphism, as is well-known in topology (the Chern characters from K-theory of tori  to $\mathbb Z_2$-graded cohomology are isomorphisms). We have therefore the following interesting consequence

Let us now consider the commutative diagram
$$\CD
K_*(C^\infty_c(G)) @>ch_G>> \PHC_*(C^\infty_c(G)))_{\mathfrak m}\\
@VVV   @VVV\\
K_{*+q}(C^\infty_c(K)) @>ch_K>> \PHC_{*+q}(C^\infty_c(K))_{\mathfrak m}\endCD$$
and 
$$\CD
K_{*+q}(C^\infty_c(K)) @>ch_K>> \PHC_{*+q}(C^\infty_c(K))_{\mathfrak m}\\
@VVV   @VVV\\
K_{*+q}(C^\infty_c(\mathbb T))^W @>ch_{\mathbb T}>> \PHC_{*+q}(C^\infty_c(\mathbb T))^W_{\mathfrak m}\endCD$$ 
Localizing $C^\infty_c(G)$ at the ideals $\mathfrak m$ , where $\mathfrak m = \{\varphi \in C^\infty_{inv}(G) | \varphi(F) = 0 \}$ for some $F \in \Omega$ ,  yields $C^\infty_c(\Omega) $, from $C^\infty_{inv}(G) \subseteq C^\infty(G)$ and so $C^\infty_c(G)_{\mathfrak m} = C^\infty_c(\Omega).$,  as  was explained in the proof of Theorem 1.12.  Then by doing similar computations to what was done in \cite{dkt1}-\cite{dkt2}, since the representations of $K$ are defined by their restrictions to maximal tori, we have the isomorphism (I), (II) and (IV). The isomorphism (III) is the main theorem 1.12.
\hfill$\Box$
\vskip .2cm

\vskip .2cm
In the next two sections we shall apply Theorem 1.12 to deduce isomorphisms of homologies of quantized algebras of functions on coadjoint orbits. In order to do this, we recall the results from \cite{diephai1} about the algebras and then compute the K-theory, periodic cyclic homology and the noncommutative 
Chern characters.

\section{Quantum half-planes}
Applying the main theorem 1.13, we compute in this section the noncommutative Chern-Connes characters for the quantum algebras of functions on the half-planes. 

\subsection{Deformation quantization}
Let us recall some results from \cite{diephai1}: 
Recall that the Lie algebra $\mathfrak g = \aff(\mathbb  R)$ of affine transformations of the real straight line is described as follows, see for example \cite{diep1}: The Lie group $\Aff(\mathbb R)$ of affine transformations:  $$x \in \mathbb R \mapsto ax + b, \mbox{ for some parameters }a, b \in \mathbb R, a \ne 0.$$ 
is known to be a two-dimensional Lie group which is isomorphic to the group of matrices
$$\Aff(\mathbb R) \cong \{\left (\begin{array}{cc} a & b \\ 0 & 1 \end{array} \right) \vert a,b \in \mathbb R , a \ne 0 \}.$$ We consider its connected component $ \Aff_0(\mathbb R)$ of the identity element given by $$G= \Aff_0(\mathbb R)= \{\left (\begin{array}{cc} a & b \\ 0 & 1 \end{array} \right) \vert a,b \in \mathbb R, a > 0 \}.$$  Its Lie algebra 
$$\mathfrak g = \aff(\mathbb R) \cong  \{\left (\begin{array}{cc} \alpha & \beta \\ 0 & 0 \end{array} \right) \vert \alpha, \beta  \in \mathbb R \}$$  admits a basis of two generators $X, Y$ with the only nonzero Lie bracket $[X,Y] = Y$, i.e. 
$$\mathfrak g = \aff(\mathbb R) \cong \{ \alpha X + \beta Y \vert [X,Y] = Y, \alpha, \beta \in \mathbb R \}.$$
The co-adjoint action of $G$ on $\mathfrak g^*$ is given (see e.g. \cite{kirillov1}) by $$\langle K(g)F, Z \rangle = \langle F, \Ad(g^{-1})Z \rangle, \forall F \in \mathfrak g^*, g \in G \mbox{ and } Z \in \mathfrak g.$$ Denote the co-adjoint orbit of $G$ in $\mathfrak g^*$, passing through $F$ by 
$$\Omega_F = K(G)F :=  \{K(g)F \vert g \in G \}.$$ Because the group $G = \Aff_0(\mathbb R)$ is exponential (see \cite{diep1}),  we have for $F \in \mathfrak g^* = \aff(\mathbb R)^*$,  
$$\Omega_F = \{ K(\exp(U)F | U \in \aff(\mathbb R) \}.$$
and hence that
$$\langle K(\exp U)F, Z \rangle = \langle F, \exp(-\ad_U)Z \rangle.$$ It is easy therefore to see that
$$K(\exp U)F = \langle F, \exp(-\ad_U)X\rangle X^*+\langle F, \exp(-\ad_U)Y\rangle Y^*.$$
For a general element $U = \alpha X + \beta Y \in \mathfrak g$, we have
$$\exp(-\ad_U) = \sum_{n=0}^\infty \frac{1}{n!}\left(\begin{array}{cc}0 & 0 \\ \beta & -\alpha \end{array}\right)^n = \left( \begin{array}{cc} 1 & 0 \\ L & e^{-\alpha} \end{array} \right),$$ where $L = \alpha + \beta + \frac{\alpha}{\beta}(1-e^\beta)$. This means that
$$K(\exp U)F = (\lambda + \mu L) X^* + (\mu e^{-\alpha})Y^*. $$ From this formula one deduces  \cite{diep1} the following description of all co-adjoint orbits of $G$ in $\mathfrak g^*$:
\begin{itemize}
\item If $\mu = 0$, each point $(x=\lambda , y =0)$ on the abscissa ordinate corresponds to a 0-dimensional co-adjoint orbit $$\Omega_\lambda = \{\lambda X^* \}, \quad \lambda \in \mathbb R .$$
\item For $\mu \ne 0$, there are two 2-dimensional co-adjoint orbits: the upper half-plane $\{(\lambda , \mu) \quad\vert\quad \lambda ,\mu\in \mathbb R , \mu > 0 \}$ corresponds to the co-adjoint orbit
\begin{equation} \Omega_{+} := \{ F = (\lambda + \mu L)X^* + (\mu e^{-\alpha})Y^* \quad \vert \quad \mu > 0 \}, \end{equation}
and the lower half-plane $\{(\lambda , \mu) \quad\vert\quad \lambda ,\mu\in \mathbb R , \mu < 0\}$ corresponds to the co-adjoint orbit
\begin{equation} \Omega_{-} := \{ F = (\lambda + \mu L)X^* + (\mu e^{-\alpha})Y^* \quad \vert \quad \mu < 0 \}. \end{equation}
\end{itemize}
We shall work henceforth on the fixed co-adjoint orbit $\Omega_+$. The case of the co-adjoint orbit $\Omega_-$ could be similarly treated. First we study the geometry of this orbit and introduce some canonical coordinates in it.
It is well-known from the orbit method \cite{kirillov1} that the Lie algebra $\mathfrak g = \aff(\mathbb R)$ is realized by the complete right-invariant Hamiltonian vector fields on co-adjoint orbits $\Omega_F \cong G_F \setminus G$ with flat (co-adjoint) action of the Lie group $G = \Aff_0(\mathbb R)$. On the orbit $\Omega_+$ we choose a fix point $F=Y^*$. It is well-known from the orbit method that we can choose an arbitrary point $F$ on $\Omega_F$. It is easy to see that the stabilizer of this (and therefore of any) point  is trivial, i.e. $G_F = \{e\}$. We identify therefore $G$ with $G_{Y^*}\setminus G$. There is a natural diffeomorphism $\Id_{\mathbb R} \times \exp(.)$ from the standard symplectic space $\mathbb R^2$ with symplectic 2-form $dp \wedge dq$ in canonical Darboux $(p,q)$-coordinates, onto the upper half-plane $\mathbb H_+ \cong \mathbb R \rtimes \mathbb R_+$ with coordinates $(p, e^q)$, which is, from the above coordinate description, also diffeomorphic to the co-adjoint orbit $\Omega_+$. We can use therefore $(p,q)$ as the standard canonical Darboux coordinates in $\Omega_{Y^*}$. There are also non-canonical Darboux coordinates $(x,y) = (p,e^q)$ on $\Omega_{Y^*}$. We show now that in these coordinates $(x,y)$, the Kirillov form looks like $\omega_{Y^*}(x,y) = \frac{1}{y}dx \wedge dy$, but in the canonical Darboux coordinates $(p,q)$, the Kirillov form is just the standard symplectic form $dp \wedge dq$. This means that there are  symplectomorphisms between the standard symplectic space $\mathbb R^2, dp \wedge dq)$, the upper half-plane $(\mathbb H_+, \frac{1}{y}dx \wedge dy)$ and the co-adjoint orbit $(\Omega_{Y^*},\omega_{Y^*})$. 
Each element $Z\in \mathfrak g$ can be considered as a linear functional $\tilde{Z}$ on co-adjoint orbits, as subsets of $\mathfrak g^*$, where $\tilde{Z}(F) :=\langle F,Z\rangle$.  It is well-known that this linear function is just the Hamiltonian function associated with the Hamiltonian vector field $\xi_Z$, which represents $Z\in \mathfrak g$ with the formula 
$$(\xi_Zf)(x) := \frac{d}{dt}f(x\exp (tZ))|_{t=0}, \forall f \in C^\infty(\Omega_+).$$ 
The Kirillov form $\omega_F$ is defined by the formula 
\begin{equation}\label{7} \omega_F(\xi_Z,\xi_T) = \langle F,[Z,T]\rangle, \forall Z,T \in \mathfrak g = \aff(\mathbb R). \end{equation} This form defines the symplectic structure and the Poisson brackets on the co-adjoint orbit $\Omega_+$. For the derivative along the direction $\xi_Z$ and the Poisson bracket we have relation $\xi_Z(f) = \{\tilde{Z},f\}, \forall f \in C^\infty(\Omega_+)$. It is well-known in differential geometry that the correspondence
$Z \mapsto \xi_Z, Z \in \mathfrak g$ defines a representation of our Lie algebra by vector fields on co-adjoint orbits. If the action of $G$ on $\Omega_+$ is flat \cite{diep1}, we have the second Lie algebra homomorphism from  strictly Hamiltonian right-invariant vector fields into the Lie algebra of smooth functions on the orbit with respect to the associated Poisson brackets.

Denote by $\psi$ the indicated symplectomorphism from $\mathbb R^2$ onto $\Omega_+$
$$(p,q) \in \mathbb R^2 \mapsto \psi(p,q):= (p,e^q) \in \Omega_+$$

It was proven in \cite{diephai1} that:
\begin{itemize}
\item Hamiltonian function $f_Z = \tilde{Z}$ in canonical coordinates $(p,q)$ of the orbit $\Omega_+$ is of the form $$\tilde{Z}\circ\psi(p,q) = \alpha p + \beta e^q, \mbox{ if  } Z = \left(\begin{array}{cc} \alpha & \beta \\ 0 & 0 \end{array} \right).$$

\item In the canonical coordinates $(p,q)$ of the orbit $\Omega_+$, the Kirillov form $\omega_{Y^*}$ is just the standard form $\omega = dp \wedge dq$.
\end{itemize}

Let us denote by $\Lambda$ the 2-tensor associated with the canonical Kirillov standard form $\omega = dp \wedge dq$ in the canonical Darboux coordinates. Recall the deformed $\star$-product of two smooth functions $u, v \in C^\infty(\Omega_\pm)$
$$u \star_h v = u.v + \sum_{r \geq 1} \frac{1}{r} \left(\frac{h}{2i}\right)^r P^r(u,v),\eqno{(\star\star)}$$ 
where $$P^r(u,v) = \Lambda^{i_1,j_1}\Lambda^{i_2,j_2}\dots \Lambda^{i_r,j_r}\partial_{i_1}\dots \partial_{i_r}u \partial_{j_1}\dots\partial_{j_r}v,$$
with the ordinary multi-index notations of the partial derivations. 
Note that in \cite{diephai1}, $(\star\star)$ was shown for the normalized Planck constant $h=1$. The situation is the same for an arbitrary nonzero value of $h$. 
It was shown that to every element $X\in \mathfrak g$ corresponds a function $\tilde{X}$ on $\mathfrak g^*$ and therefore on the K-orbits $\Omega_F$.

It was shown in \cite{diephai1}(Proposition3.1) that 
$$\frac{i}{h}\tilde X \star_h \frac{i}{h}\tilde T - \frac{i}{h}\tilde T \star_h \frac{i}{h} \tilde X = \frac{i}{h} \widetilde{[X,T]}, \forall Z,T \in \aff(\mathbb R).$$ One therefore has a representation 
$$X \longrightarrow \frac{i}{h}\tilde{X} \star_h$$
of the Lie algebra $C^\infty_c(\Omega_G)$ by the left $\star_h$-multiplication.
On the half-plane with the fixed Darboux $(q,p)$-coordinates, one fixes the Fourier transformation in $p$-coordinate 
$$\mathcal F_p(u)(\eta, q) := \frac{1}{2\pi}\int_{\mathbb R}\exp(-ip\eta)u(p,q)dp$$
and obtain that for the element $\tilde Z = \alpha p + \beta e^q$, the operator $\ell_Z$ acting on the dense subspace $L^2(\mathbb R^2, \frac{dpdq}{2\pi})^\infty$ of smooth functions by left $\star_h$-multiplication by $i\tilde Z \star_h$, i.e. $\ell_Z(u) := \frac{i}{h}\tilde Z \star_h u$. 
It was precisely computed (see Proposition 3.4 in \cite{diephai1}) that 
$$\hat{\ell}_Z(u) := \mathcal F_p \circ \ell_Z \circ \mathcal F_p^{-1}(u) = \alpha(\frac{1}{2}\partial_q - \partial_p)u.$$

It was also proven in Theorem 4.2 of \cite{diephai1} that:

The representation $\exp(\hat{\ell}_Z)$ of the group $G=\Aff_0(\mathbb R)$ is exactly the irreducible unitary representation $T_{\Omega_+}$ of $G=\Aff_0(\mathbb R)$ associated, following the orbit method construction, to the orbit $\Omega_+$, which is the upper half-plane $\mathbb H \cong \mathbb R \rtimes \mathbb R^*$, i. e.
$$(\exp(\hat{\ell}_Z)f)(y) = (T_{\Omega_+}(g)f)(y) = e^{\frac{i}{h}by}f(ay), \forall f\in L^2(\mathbb R^*, \frac{dy}{y}),$$ where $g= \exp Z = \left(\begin{array}{cc} a & b \\ 0 & 1 \end{array}\right).$

\subsection{K-groups and periodic cyclic homology and Chern-Connes characters}

Let us apply now the general notion of Chern-Connes characters to the example of the quantum algebras of functions on the coadjoint orbits of the groups of affine transformation of the real (in this section) and complex (in the next section) lines. 
Recall that the noncommutative Chern-Connes characters are some homomorphisms $ch$ from the K-groups $K_*(.)$ to the corresponding periodic cyclic homology groups $\PHC_*(.)$. 

\begin{lemma}
In the groups $\Aff(\mathbb R)$, the maximal compact subgroups is $K\cong \mathbb Z_2 = \mathbb Z /( 2 \mathbb Z )$. In its connected component of identity $\Aff_0(\mathbb R)$, the maximal compact subgroup is trivial, i.e. $K \cong \{e\}$.
\end{lemma}
{\it Proof}. Because the groups $\Aff(\mathbb R)$ is isomorphic to the semi-direct product $\mathbb R^* \ltimes \mathbb R$, one deduces that the maximal compact subgroup $K$ in $G$ is isomorphic to $\mathbb Z_2$. \hfill$\Box$

\begin{proposition} Let $\Omega_+$ be the coadjoint orbit which is the upper half-plane. Then
$$K_*(C^\infty_c(\Omega_+)) \cong \PHC_*((C^\infty_c(\Omega_+)) \cong \{e\}.$$ and therefore the noncommutative Chern-Connes characters are isomorphisms. 
\end{proposition}
{\it Proof}. Because the maximal compact subgroups of $G$ are trivial, we can conclude that the K-groups and the $\PHC_*$-groups are also trivial. \hfill$\Box$

\section{Quantum punctured complex plane}
In this section we demonstrate another application of Theorem 2.1 for the group of affine transformations of the complex line. The deformation of the coadjoint orbits of this group is in some sense more complicated than the one in the real case, see e.g. \cite{diephai2}. 

\subsection{Deformation quantization}
The group $\Aff(\mathbb C)$ is defined as $$\Aff(\mathbb C) = \{\left(\begin{array}{cc} a & b \\ 0 & 1 \end{array}\right) \vert \quad a, b \in \mathbb C, a \ne 0 \}.$$ It is isomorphic to the semi-direct product of the complex line $\mathbb C$ and the punctured complex line $\mathbb C^* = \mathbb C \setminus (0)$. The group is connected but not simply connected. and the exponent map $$\exp : \mathbb C \to \mathbb C^*; z \mapsto e^z$$ gives rise to the universal covering 
$$\widetilde{\Aff(\mathbb C)} \cong \mathbb C \rtimes \mathbb C \cong \{(z,w) \vert z,w \in \mathbb C\}$$ of $\Aff(\mathbb C)$ with multiplication
$$(z,w)(z',w') := (z+z',w +e^zw').$$ 
As a real Lie group, it is 4-dimensional and we denote its Lie algebra by $\aff(\mathbb C) = \Lie\Aff(\mathbb C)$ The dual space $\mathfrak g^*$ of $\mathfrak g = \Lie\Aff(\mathbb C)$ can be identify with $\mathbb R^4$ with coordinates $(\alpha, \beta, \gamma, \delta)$,  see \cite{diep1}. The coadjoint orbits of $\widetilde{\Aff(\mathbb C)}$ in $\mathfrak g^*$ passing though a point $F= \alpha X_1^* + \beta X_2^* + \gamma Y_1^* + \delta Y_2^*$ is denote by $\Omega_F$, where $X_1^*, X_2^*, Y_1^*, Y_2^*$ form the basis of $\mathfrak g^*$ dual to the basis $X_1, X_2, Y_1, Y_2$ of $\mathfrak g$ with the brackets
$$[X_1,Y_1]= Y_1, [X_1,Y_2]= Y_2, [X_2,Y_1]= Y_2, [X_2,Y_2]=-Y_1.$$ Then 
\begin{itemize}
\item Each point $(\alpha,0,0,\delta)$ is a 0-dimensional coadjoint orbit, denoted $\Omega_{(\alpha,0,0,\delta)}$,
\item The open set $\beta^2 + \gamma^2 \ne 0$ is the single 4-dimensional coadjoint orbit $\Omega \approx \mathbb C \times \mathbb C^*$, the punctured complex plane. 
\end{itemize}
Note that the orbit $\Omega$ is not simply connected and there is no diffeomorphism from some symplectic vector space onto it. In \cite{diephai2}, some system of diffeomorphisms was constructed. Let us recall them. Consider 
$$\mathbb H_k := \{w=q_1+iq_2 \in \mathbb C \vert -\infty < q_1 < +\infty, 2k\pi < q_2 < 2k\pi < q_2 < (2k+1)\pi \},$$ for each $k=0,\pm 1,\dots$. Let $\mathbb C_k := \mathbb C \setminus L$,  where $L$ is the positive real line 
$$L = \{\rho e^{i\varphi} \in \mathbb C \vert 0 < \rho < \infty, \varphi = 0 \}.$$ There is a natural map $$ \mathbb C \times \mathbb C \to \Omega \cong \mathbb C \times \mathbb C^*; (z,w) \mapsto (z, e^w),$$ whose restriction gives a diffeomorphism
$$\varphi_k : \mathbb C \times \mathbb H_k \to \mathbb C \times \mathbb C^*.$$
On $\mathbb C \times \mathbb H_k$ we have the natural symplectic form
$$\omega_0 := \frac{1}{2}[dz \wedge dw + d\bar{z} \wedge d\bar{w}],$$ induced from the standard symplectic form on $\mathbb C^2$ with coordinates $(z,w)$. The corresponding symplectic form matrix is 
$$\wedge = \left[\begin{array}{cccc} 0 & -1 &  0 & 0\\ 1 & 0 & 0 & 0 \\ 0 & 0 & 0 & 1 \\  0 & 0 &  -1 & 0 \end{array}\right] \mbox{ and } \wedge^{-1} = \left[\begin{array}{cccc} 0 & 1 &  0 & 0\\ -1 & 0 & 0 & 0 \\ 0 & 0 & 0 & -1 \\  0 & 0 &  1 & 0 \end{array}\right] .$$ The corresponding Poisson brackets of functions $f,g\in C^\infty(\Omega)$ is 
$$\{f,g\} = \wedge^{ij} \frac{\partial f}{\partial x^i}\frac{\partial g}{\partial x^j} = \frac{\partial f}{\partial p_1}\frac{\partial g}{\partial q_1} - \frac{\partial f}{\partial q_1}\frac{\partial g}{\partial p_1} - \frac{\partial f}{\partial p_2}\frac{\partial g}{\partial q_2} + \frac{\partial f}{\partial q_2}\frac{\partial g}{\partial p_2}.$$
For an arbitrary element $A \in  \aff(\mathbb C)$ it was computed in \cite{diephai2}(Proposition 2.4) that:
\begin{itemize}
\item the corresponding function on $\Omega$ is 
$$\tilde{A}\circ \varphi_k(z,w) = \frac{1}{2}[\alpha z + \beta e^w + \bar{\alpha}\bar{z} + \bar{\beta} e^{\bar{w}}].$$
\item In the local coordinates $(z,w)$ of the orbit $\Omega$, the Kirillov form $\Omega$  coincides with the standard form
$$\omega_0 := \frac{1}{2}[dz \wedge dw + d\bar{z} \wedge d\bar{w}].$$
\end{itemize}
It was also proven in \cite{diephai2}(Proposition 3.1) that for all $A, B \in \aff(\mathbb C)$, the Moyal $\star_h$ product satisfies the relation
$$\frac{i}{h}\tilde{A} \star_h \frac{i}{h}\tilde{B} - \frac{i}{h}\tilde{B} \star_h \frac{i}{h}\tilde{A} = \frac{i}{h}[A,B].$$ This means that we have some representation  $\ell^{(k)}_A: A \mapsto \frac{i}{h}\tilde{A} \star_h$ of Lie algebra $\aff(\mathbb C)$ on the space $C^\infty(\Omega)$. 
Denote by $\mathcal F_z$ the Fourier transformation
$$\mathcal F_z(f)(\xi,w) := \frac{1}{2\pi}\int_{\mathbb R^2} \exp(-i\Re(\xi\bar{z}))f(z,w)dp_1dp_2$$ and the inverse Fourier transformation
$$\mathcal F_z^{-1}(\tilde{f})(z,w) := \frac{1}{2\pi}\int_{\mathbb R^2} \exp(i\Re(\xi\bar{z}))\tilde{f}(\xi,w)dp_1dp_2$$
By the same computation as in \cite{diephai2} we have
for each $A = \left(\begin{array}{cc} \alpha & \beta \\ 0 & 0 \end{array}\right) \in \aff(\mathbb C)$ and for each compactly supported smooth function $f\in C^\infty_c(\mathbb C \times \mathbb H_k)$, 
$$\hat{\ell}_Af = \mathcal F_z \circ \ell^{(k)}_A \circ \mathcal F^{-1}_z(f) = [\alpha(\frac{1}{2}\partial_w - \partial_{\bar{\xi}}) + \bar{\alpha}(\frac{1}{2}\partial_{\bar{w}} - \partial_{\xi}) + \frac{i}{2}(\beta e^{w - \frac{1}{2}\bar{\xi}} + \bar{\beta} e^{\bar{w} - \frac{1}{2\xi}})f].$$
It was shown \cite{diephai2}(Theorem 4.2) that the representation $\exp(\hat{\ell}^{(k)}_A)$ of the universal covering 
group $\widetilde{\Aff(\mathbb C)}$  coincides with the irreducible unitary representation $T_\theta$ of $\widetilde{\Aff(\mathbb C)}$ associated with $\Omega$ by the orbit method, i.e. 
$$\exp(\hat{\ell}^{(k)}_A)f(x) = [T_\theta(\exp A)f](x),$$ realizing the space $L^2(\mathbb R \times \mathbb S^1)$ and acting as 
$$[T_\theta(z,w)f](x) = \exp(\frac{i}{h}\Re(wx) + 2\pi\theta[\frac{\Im(x +z)}{2\pi}]))f(x \oplus z),$$
where $(z,w)\in \widetilde{\Aff(\mathbb C)}$, $x\in \mathbb R \times \mathbb S^1 = \mathbb C \setminus (0)$, $f \in L^2(\mathbb R \times \mathbb S^1),$ $x \oplus z := \Re(x+z) + 2\pi i \{\frac{x+z}{2\pi}\},$ $[a]$ is the integral part of $a$ and $\{a\}$ is the decimal part of $a$.

\subsection{K-groups and periodic cyclic homology}

\begin{lemma}
The maximal compact subgroup
$K = \{ \left[\begin{array}{cc} e^{i\varphi} & 0 \\ 0 & 1 \end{array}\right] \}$
 of $\Aff(\mathbb C)$ is isomorphic to $\mathbb S^1$. 
\end{lemma}
{\it Proof} is easy and is omitted. \hfill$\Box$ 

From this one deduces the following results.
\begin{proposition} Let $\Omega$ b the coadjoint orbit of $\Aff(\mathbb C)$, which is the punctured complex plane. Then,
$$K_0(C^\infty_c(\Omega)) \cong \mathbb Z \mbox{ and } K_1(C^\infty_c(\Omega)) \cong \{0\}.$$
$$\PHC_0(C^\infty_c(\Omega)) \cong \mathbb Z \mbox{ and } \PHC_1(C^\infty_c(\Omega)) \cong \{0\}.$$
\end{proposition}
{\it Proof}. Because of Theorem 1.13 and Lemma 4.1, $K_*(C^\infty_c(\Omega))$ are the same as $K_*(C^\infty(\mathbb S^1)) \cong K^*(\mathbb S^1)$, and $\PHC_*(C^\infty_c(\Omega))$ are isomorphic to $\PHC_*(C^\infty_c(\mathbb S^1))$. The assertions become clear.   \hfill$\Box$

\subsection{Chern-Connes characters}

\begin{proposition} Let $\Omega$ be the coadjoint orbit of $\Aff(\mathbb C)$, which is the punctured complex plane. 
Then, the Chern-Connes character $$ch: K_*(C^\infty_c(\Omega)) \to \PHC_*(C^\infty_c(\Omega))$$ is an isomorphism. 
\end{proposition}
{\it Proof}. 
From Theorem 2.1, we have the commutative diagram:
$$\CD
K_*(C^\infty_c(\Omega)) @>ch_\Omega>> \PHC_*(C^\infty_c(\Omega)))\\
@VIVV   @VVIIIV\\
K_{*+q}(C^\infty_c(\Omega_K)) @>ch_{\Omega_K}>> \PHC_{*+q}(C^\infty_c(\Omega_K))\\
@VIIVV   @VV{IV}V\\
K_{*+q}(C^\infty_c(\mathbb T))^W @>ch_{\mathbb T}>> \PHC_{*+q}(C^\infty_c(\mathbb T))^W_{\mathfrak m}
\endCD$$
It was shown \cite{dkt1}-\cite{dkt2} that the Chern-Connes characters are reduced to the classical Chern characters of commutative tori. For the tori, $ch_{\mathbb T}$ is an isomorphism modulo torsions, and therefore $ch_{\Omega_K}$ and $ch_\Omega$ are also isomorphisms modulo torsions. In our case of the quantum punctured complex plane, the groups are either 0 or $\mathbb Z$. Hence, Chern-Connes character
$$ch_\Omega : K_*(C^\infty_c(\Omega)) \to \PHC_*(C^\infty_c(\Omega))$$
 is  an isomorphism.\hfill$\Box$

\section{Quantum coadjoin orbits of $\SL(2,\mathbb R)$}

\subsection{Orbit method for $\SL(2,\mathbb R)$}

The group $\SL(2,\mathbb R)$ is well-known as the special linear groups the Lie algebra of which is generated by 3 Cartan basis elements
$$H = \left[\begin{array}{cc} 1 & 0\\ 0 & -1\end{array}\right],
X = \left[\begin{array}{cc} 0 & 1\\ 1 & 0\end{array}\right], \mbox{ and }
Y = \left[\begin{array}{cc} 0 & 1\\ -1 & 0\end{array}\right], $$
subject to the relation
$$[H,X] = 2Y,\quad [H,Y] =2X,\quad [X,Y] = -2H.$$ Let us denote by $H^*, X^*, Y^*$ the dual basis of $\mathfrak g^*$. The Killing form on the Lie algebra $\sl(2,\mathbb R)$ of $SL(2,\mathbb R)$ is defined by the ordinary trace
$\langle X,Y\rangle := \frac{1}{4}\tr(\ad X\circ \ad Y)$.
It is easy to see that by the real isomorphism (complex anti-isomorphism) $\mathfrak g \to \mathfrak g^*$, defined by $Y \mapsto \langle .,Y\rangle$ we can identify $\mathfrak g$ with $\Ad$ action of $G$, with $\mathfrak g^*$ with coadjoint action $K(g) = \Ad^*(g^{-1})$ of $G$. This also mean that $\widehat{\Ad(g)X} = \Ad^*(g^{-1})\hat{X}.$ The problem is therefore reduced to compute the conjugacy classes in matrices. It is easy to show that each matrix in $\sl(2,\mathbb R)$ can be made conjugae to one of the following standard matrices
$$\left[\begin{array}{cc} 0 & \lambda\\ -\lambda & 0\end{array}\right],
\left[\begin{array}{cc} 0 & 0\\ 1 & 0\end{array}\right],
\left[\begin{array}{cc} -\lambda & 0\\ 0 &\lambda \end{array}\right],
\left[\begin{array}{cc} 0 & 0\\ 0 & 0\end{array}\right].$$
and their $W_2$ conjugates
$$
\left[\begin{array}{cc} 0 & -\lambda\\ \lambda & 0\end{array}\right],
\left[\begin{array}{cc} 0 & 1\\ 0 & 0\end{array}\right],
\left[\begin{array}{cc} \lambda & 0\\ 0 &-\lambda \end{array}\right],
$$ 
where $W_2=\mathbb Z_2$ is the Weyl group of $\SL(2,\mathbb R)$ acting on the Cartan subalgebra of $\sl(2,\mathbb R)$ by changing the basis.

It is easy to check out, see e.g. \cite{doduchanh} the list of all coadjoint orbits of $\SL(2,\mathbb R)$:
\begin{itemize}
\item Elliptic hyperboloid:
$$\Omega^1_\lambda = \{ 2xX^* + 2hH^* + 2yY^* | x^2 + h^2 - y^2 = \lambda^2, \lambda \ne 0 \},$$
\item Upper half-cones:
$$\Omega^2_+ =\{ 2xX^* + 2hH^* -2yY^* | x^2 +h^2 - y^2 = 0, y > 0 \},$$

Lower half-cones: 
$$\Omega^2_+ =\{ 2xX^* + 2hH^* -2yY^* | x^2 +h^2 - y^2 = 0, y < 0 \},$$

One-point:
$$\Omega^2_0 = \{0\},$$

\item
Upper elliptic half-hyperboloid:
$$\Omega^3_{\lambda,+} = \{ 2xX^* + 2hH^* - 2yY^* | x^2 + h^2 - y^2 = -\lambda^2, y > 0 \},$$
Lower elliptic half-hyperboloid:
$$\Omega^3_{\lambda,-} = \{ 2xX^* + 2hH^* - 2yY^* | x^2 + h^2 - y^2 = -\lambda^2, y < 0 \}.$$
\end{itemize}

In order to obtain unitary representations associated with coadjoint orbit, following the orbit method, we should choose for every orbit a complex subalgebra in the complex hull of the Lie algebra $\mathfrak g = \sl(2,\mathbb R)$ satisfying conditions of a polarization, see e. g. \cite{diep1}, \cite{kirillov1}.

\begin{itemize}
\item For the orbit $\Omega^1_\lambda$, fix a point $F=2\lambda H^*$. The complex subalgebra $\mathfrak p := \mathbb C H + \mathbb C (X+Y)$ is a polarization with $\chi_F = \exp(2\pi i \langle F,.\rangle)$ and the polarizing subgroup $H = H_0 \cup \varepsilon H_0$ and the extended character $\widetilde{\chi}_F(\varepsilon) = \pm 1$.
\item For the orbit $\Omega^2_{\lambda,\pm}$, fix a point $F=X^* - Y^*$. The complex subalgebra $\mathfrak p := \mathbb C H + \mathbb C (X+Y)$ is a polarization with $\chi_F = \exp(2\pi i \langle F,.\rangle)$ and the polarizing subgroup $H = H_0 \cup \varepsilon H_0$ and the extended character $\widetilde{\chi}_F(\varepsilon) = \pm 1$.
\item For the orbit $\Omega^3_{\lambda,\pm}$, fix a point $F=2H^*$. The complex subalgebra $\mathfrak p := \mathbb C H + \mathbb C (X+iH)$ is a polarization with $\chi_F = \exp(2\pi i \langle F,.\rangle)$ and the polarizing subgroup $H = H_0\rtimes \SO(2,\mathbb R)$ and the extended character 
$$\widetilde{\chi}_F(\left[\begin{array}{cc} \cos a & \sin a\\ -\sin a & \cos a\end{array}\right]) = \exp(-4\pi i\lambda a),$$ if and only if $\lambda = \frac{k}{8}$ for a integer $k\in \mathbb Z$: The stabilizer $G_F$ is in this case $\SO(2,\mathbb R)$.
\end{itemize}

Following the orbit method, see \cite{diep1},\cite{kirillov1}, after choosing polarizations, the quantum induced bundles are constructed and the representations are obtained as the natural actions of the group on quantum bundles. Remark that the polarizing subalgebras contain all the ``p'' ordinates and the quotient spaces are realized on ``q''ordinates in the Darboux coordinates associated with polarizations. The action of $\star$-products on smooth sections gives us the quantum algebra structure on the spaces of smooth sections of the quantum bundles. The $\star$-product are constructed following the method of Fedosov quantization.

\subsection{Deformation quantization on coadjoint orbits}
Each coadjoint orbit of $\SL(2,\mathbb R)$ is a symplectic manifold:
\begin{itemize}
\item For the case of $\Omega^1_\lambda$, we can choose local coordinates such that the map
$$\psi : (p,q) \in \mathbb R^2 \mapsto (x,h,y) \in \Omega^1_\lambda$$ such that
$$\left\{\begin{array}{ccl} 
x & = & p\cos q - \lambda \sin q,\\
h & = & p\sin q + \lambda \cos q,\\
y & = & p 
\end{array}\right.$$ is a local symplectomorphism.
An element $A\in \mathfrak g$ can be considered as a Hamiltonian $\tilde{A}\in C^\infty(\Omega^1_\lambda)$ on the coadjoint orbit as a symplectic manifold. In the introduced Darboux coordinates $(p,q)$ we have \cite{doduchanh}
$$\tilde{A}(F) = \langle F,A\rangle = (2a_1 \cos q + 2b_1 \sin q - 2c_1)p + (-2a_1 \sin q + 2b_1 \cos q)\lambda.$$   
\item
For the orbits $\Omega^3_{\lambda,\pm}$ we can choose local symplectomorphism
$$\psi : \mathbb C \times \mathbb C \to (\Omega^3_{\lambda, \pm})_{\mathbb C}$$with the same $\tilde{A}$ and
\item For the orbits $\Omega^2_{\pm}$ we can choose the same as in case of $\Omega^1_{\lambda=0}$ or $\Omega^3_{\lambda = 0,\pm}$ also with the same $\tilde{A}$.
\end{itemize}

The deformation quantization operators
$$\ell_f : C^\infty_c(\Omega) \to C^\infty(\Omega), $$ defined by 
$$g\mapsto f \star g$$ can be exactly computed as some differential operators
For $s = q - \frac{x}{2}$, the Fourier transform in variable $p$ is
$$\hat{\ell}_A = \mathcal F_p \circ \ell_A \circ \mathcal F_p^{-1} (f) = (a_1\cos s + b_1\sin s - c_1) \partial_s + (-a_1\sin s + b_1\cos s)(2\lambda i + 1).$$
The corresponding unitary representation of the universal covering $\widetilde{\SL(2,\mathbb R)}$ of $\SL(2,\mathbb R)$ is obtained as the unitary operator-valued solution of the Cauchy problem
$$\left\{\begin{array}{rcl}\partial_tU(\exp(A),t) &=& \ell_A U(\exp(A),t),\\
U(\exp(A),t)|_{t=0} &=& Id.\end{array}\right.$$
They are exactly the principal series, discrete series irreducible representations of $\SL(2,\mathbb R)$. The remaining four irreducible unitary representations $\pi_0^{\pm}$, $\pi_1^{\pm}$ are obtained by taking the limits  of the discrete series $\pi_s$,of the parameters $0<s<1$, $s \to 0$ or $ s \to 1$ and decomposing the limits in sums of two irreducible components. These representations ``correspond'' to the nilpotent orbits, half-cones or the ones passing through $F = \left[\begin{array}{cc} 0 & 0 \\ 1 & 0\end{array}\right]$ or $F = \left[\begin{array}{cc} 0 & 1 \\ 0 & 0\end{array}\right]$

Together with $\star$-product we have the quantized coadjoint orbits as explained in \cite{doduchanh}. We compute in the rest of this paper the corresponding K-theory, HP-theory and noncommutative Chern-Connes characters.

\subsection{Quantum elliptic hyperboloids}
In this subsection we study the quantum elliptic hyperboloids $(C^\infty_c(\Omega^1_{\lambda,\pm}),\star_\hbar))$ and compute the K-theory, HP-theory and the Chern character.
\begin{proposition}
$$K_*(C^\infty_c(\Omega^1_{\lambda,\pm}),\star_\hbar)\cong K^*(\mathbb S^1),$$
$$\PHC_*(C^\infty_c(\Omega^1_{\lambda,\pm}),\star_\hbar)\cong \PHC^*(\mathbb S^1).$$
\end{proposition}
{\it Proof.} 
The maximal compact subgroup of $\SL(2,\mathbb R)$ is $K=\SO(2,\mathbb R)$ and $\Omega^1_\lambda|_K$ is $\SO(1,\mathbb R) \setminus \SO(2,\mathbb R) \cong \SO(2,\mathbb R)\cong \mathbb S^1$ so, from  Theorem 2.1, we have
$$K_*(C^\infty_c(\Omega^1_{\lambda,\pm}),\star_\hbar)\cong K_{*+2}(C^\infty_c(\Omega^1_{\lambda,\pm})|_K),\star_\hbar)\cong K^*(\mathbb S^1),$$
$$\PHC_*(C^\infty_c(\Omega^1_{\lambda,\pm}),\star_\hbar)\cong \PHC_{*+2}(C^\infty_c(\Omega^1_{\lambda,\pm}|_K),\star_\hbar)\cong \PHC^*(\mathbb S^1).$$
\hfill$\Box$

\begin{proposition}
The Chern-Connes character for the elliptic hyperboloids  $$ch: K_*(C^\infty_c(\Omega^1_{\lambda,\pm}),\star_\hbar) \to \PHC_*(C^\infty_c(\Omega^1_{\lambda,\pm}),\star_\hbar)$$ is an isomorphism.
\end{proposition}
{\it Proof.} In this case the Chern character of topological torus $\mathbb T \cong \mathbb S^1$ is an isomorphism and we have $$\CD ch: K_*(C^\infty_c(\Omega^1_{\lambda,\pm})) \cong K^*(\mathbb S^1) @>\cong >> H^*_{DR}(\mathbb S^1) \cong \PHC_*(\mathbb S^1)\cong \PHC_*(\Omega^1_{\lambda,\pm})\endCD$$ is also an isomorphism.
\hfill$\Box$
\vskip 0.5cm

It is easy to see that for the other quantum coadjoint orbits we have the same results.

\subsection{Quantum  elliptic two-fold hyperboloids}
Let us now study the quantum elliptic two-fold hyperboloids $(C^\infty_c(\Omega^3_{\lambda,\pm}),\star_\hbar))$ and compute the K-theory, HP-theory and the Chern character.
\begin{proposition}
$$K_*(C^\infty_c(\Omega^3_{\lambda,\pm}),\star_\hbar)\cong K^*(\mathbb S^1),$$
$$\PHC_*(C^\infty_c(\Omega^3_{\lambda,\pm}),\star_\hbar)\cong \PHC^*(\mathbb S^1)$$
\end{proposition}
{\it Proof.} 
From our notation, if we fix a point on the orbit $\Omega^3_{\lambda,\pm}$, then $(\Omega^3_{\lambda,\pm})|_K $ is isomorphic to $\SO(1)\setminus \SO(2)\approx \mathbb S^1$ and we have again
$$K_*(C^\infty_c(\Omega^3_{\lambda,\pm}),\star_\hbar)\cong K_{*+2}(C^\infty_c(\Omega^3_{\lambda,\pm}|_K),\star_\hbar)\cong K^*(\mathbb S^1).$$
$$\PHC_*(C^\infty_c(\Omega^3_{\lambda,\pm}),\star_\hbar)\cong \PHC_{*+2}(C^\infty_c(\Omega^3_{\lambda,\pm}|_K),\star_\hbar)\cong \PHC^*(\mathbb S^1).$$

\hfill$\Box$
\begin{proposition}
The Chern-Connes character for the elliptic two-fold hyperboloids  $$ch: K_*(C^\infty_c(\Omega^3_{\lambda,\pm}),\star_\hbar) \to \PHC_*(C^\infty_c(\Omega^3_{\lambda,\pm}),\star_\hbar)$$ is an isomorphism.
\end{proposition}
{\it Proof.} We have again that the Chern character of topological torus $\mathbb T \cong \mathbb S^1$ is an isomorphism and we have $$\CD ch: K_*(C^\infty_c(\Omega^3_{\lambda,\pm})) \cong K^*(\mathbb S^1) @>\cong >> H^*_{DR}(\mathbb S^1) \cong \PHC_*(\mathbb S^1)\cong \PHC_*(\Omega^3_{\lambda,\pm}).\endCD$$
\hfill$\Box$

\subsection{Quantum elliptic half-cones}
We now study the quantum elliptic half-cones $(C^\infty_c(\Omega^2_{\lambda,\pm}),\star_\hbar))$ and compute the K-theory, HP-theory and the Chern character.
\begin{proposition}
$$K_*(C^\infty_c(\Omega^2_{\lambda,\pm}),\star_\hbar)\cong K^*(\mathbb S^1),$$
$$\PHC_*(C^\infty_c(\Omega^2_{\lambda,\pm}),\star_\hbar)\cong \PHC^*(\mathbb S^1)$$ is an isomorphism.
\end{proposition}
{\it Proof.} Fix a point on the orbit $\Omega^2_{\lambda,\pm}$ which is a half-cone and by a direct computation we have $(\Omega^2_{\lambda,\pm})|_K \cong \SO(1)\setminus \SO(2)\approx \mathbb S^1$. \hfill$\Box$

\begin{proposition}
The Chern-Connes character for the elliptic half-cones  $$ch: K_*(C^\infty_c(\Omega^2_{\lambda,\pm}),\star_\hbar) \to \PHC_*(C^\infty_c(\Omega^2_{\lambda,\pm}),\star_\hbar)$$ is an isomorphism.
\end{proposition}
{\it Proof.} In this case the Chern character of topological torus $\mathbb T \cong \mathbb S^1$ is an isomorphism and we have $$\CD ch: K_*(C^\infty_c(\Omega^2_{\lambda,\pm})) \cong K^*(\mathbb S^1) @>\cong >> H^*_{DR}(\mathbb S^1) \cong \PHC_*(\mathbb S^1)\cong \PHC_*(\Omega^2_{\lambda,\pm}).\endCD$$
\hfill$\Box$

\section*{Acknowledgments}
This work was completed during the time the first author visited Abdus Salam ICTP in 2001 and 2003. The first author expresses his sincere thanks for the hospitality and for the excellent conditions provided.

{\sc Institute of Mathematics, National Centre for Science and Technology of Vietnam, P. O. Box 631, Bo Ho 10.000, Hanoi Vietnam}

Email: {\tt dndiep@thevinh.ncst.ac.vn}

\vskip .5cm
{\sc Abdus Salam International Centre for Theoretical Physics, Strada Costiera 11, 34014, Trieste, Italy}

Email: {\tt kuku@ictp.trieste.it}
\end{document}